# LEARNING MIXTURES OF SEPARATED NONSPHERICAL GAUSSIANS

By Sanjeev Arora[1]  and  Ravi Kannan[2]

*Princeton University and Yale University*

Mixtures of Gaussian (or normal) distributions arise in a variety of application areas. Many heuristics have been proposed for the task of finding the component Gaussians given samples from the mixture, such as the *EM algorithm*, a local-search heuristic from Dempster, Laird and Rubin [*J. Roy. Statist. Soc. Ser. B* **39** (1977) 1–38]. These do not provably run in polynomial time.

We present the first algorithm that provably learns the component Gaussians in time that is polynomial in the dimension. The Gaussians may have arbitrary shape, but they must satisfy a "separation condition" which places a lower bound on the distance between the centers of any two component Gaussians. The mathematical results at the heart of our proof are "distance concentration" results—proved using isoperimetric inequalities—which establish bounds on the probability distribution of the distance between a pair of points generated according to the mixture.

We also formalize the more general problem of max-likelihood fit of a Gaussian mixture to unstructured data.

**1. Introduction.** Finite mixture models are ubiquitous in a host of areas that use statistical techniques, including artificial intelligence (AI), computer vision, medical imaging, psychology and geology (see [15, 23]). A mixture of distributions $\mathcal{D}_1, \mathcal{D}_2, \ldots$ with mixing weights $w_1, w_2, w_3, \ldots$ (where $\sum_i w_i = 1$) is the distribution in which a sample is produced by first picking a component distribution—the $i$th one is picked with probability $w_i$—and then producing a sample from that distribution. In many applications the component distributions are multivariate Gaussians.

Received September 2003; revised December 2003.

[1]Supported by NSF Grants CCR-00-98180 and CCR-02-05594, an Alfred Sloan Fellowship and a David and Lucille Packard Fellowship.

[2]Supported in part by NSF Grant CCR-98-20850.

Given samples from the mixture distribution, how can one learn (i.e., reconstruct) the component distributions and their mixing weights? The most popular method is probably the EM algorithm of Dempster, Laird and Rubin [7]. EM is a local search heuristic that tries to converge to a *maximum-likelihood* description of the data by trying to cluster sample points according to which Gaussian they came from. Though moderately successful in practice, it often fails to converge or gets stuck in local optima. Much research has gone into fixing these problems, but has not yet resulted in an algorithm that provably runs in polynomial time. A second known technique is called *projection pursuit* in statistics [12]. In this, one projects the samples into a random low-dimensional space and then, in the projected space, tries to do the clustering (exhaustively) exploiting the low dimensionality.

We note that some combinatorial problems seemingly related to learning Gaussian mixtures are NP-hard. For instance, Megiddo [18] shows that it is NP-hard to decide, given a set of points in $\Re^n$, whether the points can be covered by two unit spheres. This problem seems related to learning a mixture of two spherical Gaussians.

Nevertheless, one may hope that when the data is generated from the mixture of Gaussians (as opposed to being unstructured as in Megiddo's result) then the algorithm could use this structure in the data. Recently, Dasgupta [5] took an important step in this direction by showing how a mixture of $k$ identical Gaussians could be learned in polynomial time provided the Gaussians are "spherelike" (their probability mass is concentrated in a thin spherical shell) and their centers are "well-separated." (Such separation conditions correspond to a nondegeneracy assumption: if the mixture contains two identical Gaussians whose centers are arbitrarily close, then they cannot be distinguished even in principle.)

Though Dasgupta's algorithm is a good first step, it leaves open the question whether one can design algorithms that require weaker assumptions on the Gaussians. This is the topic of the current paper: our algorithms make no assumption about the shape of the Gaussians but they require the Gaussians to be "well-separated." Even for the special case of spherical Gaussians, our result improves Dasgupta's (and a result of Dasgupta and Schulman [6] that is independent of our work). We describe our results in more detail in Section 2.3 and compare them to other work.

We also define a more general problem of Gaussian fitting, whereby we make no assumptions about the data and have to fit the mixture of $k$ Gaussians that maximizes the log-likelihood it assigns to the dataset (see Section 2.1). We use techniques developed in the context of approximation algorithms to design algorithms for one of the problems (see Section 4). The exact problem is NP-hard.

**2. Definitions and overview.** The univariate distribution $N(\mu, \sigma)$ on $\Re$ has the density function $f(x) = (\sqrt{2\pi}\sigma)^{-1} \exp(-\frac{(x-\mu)^2}{2\sigma^2})$. It satisfies $E[(x - \mu)^2] = \sigma^2$. The analogous distribution in $\Re^n$ is the *axis-aligned* Gaussian $N(\bar{\mu}, \bar{\sigma})$, where $\bar{\mu}, \bar{\sigma} \in \Re^n$ and the density function is the product distribution of $N(\mu_1, \sigma_1), N(\mu_2, \sigma_2), \ldots, N(\mu_n, \sigma_n)$. A random sample $(x_1, x_2, \ldots, x_n)$ satisfies $E[\sum_i (x_i - \mu_i)^2] = \sum_i \sigma_i^2$. (Similarly, $E[\sum_i (x_i - \mu_i)^2 / \sigma_i^2] = n$.)

A general Gaussian in $\Re^n$ is obtained from an axis-aligned Gaussian by applying an arbitrary rotation. Specifically, its probability density function has the form

$$(1) \qquad F_{Q,p}(x) = \frac{1}{(2\pi)^{n/2} \prod_i \sqrt{\lambda_i(Q)}} \exp(-(x-p)^T Q^{-1} (x-p)/2),$$

where $Q$ is an $n \times n$ positive definite matrix with eigenvalues $\lambda_1(Q), \ldots, \lambda_n(Q) > 0$, and $p \in \Re^n$ is the *center*. Since $Q$ can be rewritten as $R^{-1} \times \text{diag}(\lambda_i(Q))R$, where $R$ is a rotation, the quantities $\lambda_i(Q)$ play the same role as the variances $\sigma_i^2$ in the axis-aligned case. From our earlier discussion, $E[(x-p)^T \times (x-p)]$ is $\sum_i \lambda_i(Q)$ and $E[(x-p)^T Q^{-1}(x-p)] = \int_x F_{Q,p}(x)(x-p)^T Q^{-1}(x-p) = n$.

For any finite sample of points in $\Re^n$ we can try to fit a Gaussian by computing their *variance–covariance matrix*. Let $x_1, x_2, \ldots, x_N$ be $N$ points in $\Re^n$ in general position (i.e., we assume that they do not all lie on a hyperplane). Let $X$ be the $n \times N$ matrix whose columns are the vectors $x_1 - q, x_2 - q, \ldots, x_N - q$, where $q = \frac{1}{N}(x_1 + x_2 + \cdots + x_N)$ is the sample mean. Then the variance–covariance matrix $A = \frac{1}{N} X X^T$; note that it is positive definite by definition.

This fit may, of course, be poor for an arbitrary point set. However, for every $\varepsilon > 0$, there is a constant $c_\varepsilon > 0$ such that if $N \geq c_\varepsilon n \log n$ and the $N$ points are independent, identically distributed samples from a Gaussian $F_{G,p}$, then with probability at least 0.99, $F_{A,q}$ is a $(1+\varepsilon)$-fit to $F_{G,p}$ in every direction [3, 22] in the sense that $|q - p|^2 \leq \varepsilon \sum_i \lambda_i(Q)$ and $|Gv|(1 - \varepsilon) \leq |Av| \leq (1+\varepsilon)|Gv|$ for every unit length vector $v$. (The proof of this is highly nontrivial; but a weaker statement, when the hypothesis is strengthened to $N \geq c_\varepsilon n^2$, is easier to prove.)

*Spherical and spherelike Gaussians.* In an axis-aligned Gaussian with center at the origin and with variances $\sigma_1^2, \ldots, \sigma_n^2$, the quantity $\sum_i x_i^2 / \sigma_i^2$ is the sum of $n$ independent identical random variables from $N(0, 1)$ so this sum is tightly concentrated about its mean $n$. In a *spherical* Gaussian, all $\sigma_i$'s are the same, so even $\sum_i x_i^2$ is tightly concentrated. (These observations go back to Borel.) More generally, $E[\sum_i x_i^2] = \sum_i \sigma_i^2$. If the $\sigma_i$'s are not "too different," then distance-concentration results (similar to Lemma 5 below) show that that almost all of the probability mass is concentrated in a thin spherical shell of radius about $(\sum_i \sigma_i^2)^{1/2}$; such Gaussians may be thought of as *spherelike*. Roughly speaking, if radius/$\sigma_{\max} = \Omega(\log n)$, then the Gaussian is spherelike. Known algorithms (such as [5]) work only for such spherelike Gaussians. By contrast, here, we wish to allow Gaussians of all shapes.

2.1. *Max-likelihood learning.* Now we formalize the learning problems. Consider a mixture of Gaussians $(w_1, F_1, w_2, F_2, \ldots, w_m, F_m)$ in $\Re^n$, where the $w_i$'s are the mixing weights. With any point $x \in \Re^n$, one can associate $m$ numbers $(F_i(x))_{i=1,\ldots,m}$ corresponding to the probabilities assigned to it by the various Gaussians according to (1). For any sample $S \subseteq \Re^n$ this imposes a natural partition into $m$ blocks: each point $x \in S$ is labeled with an integer $l(x) \in \{1, \ldots, m\}$ indicating the distribution that assigns the highest probability to $x$. (Ties are broken arbitrarily.) The *likelihood* of the sample is

$$\prod_{x \in S} F_{l(x)}(x).$$

It is customary to work with the logarithm of this likelihood, called the *log-likelihood.*

Thus one may mean several things when talking about "learning mixtures of Gaussians" [21]. The following is the most general notion.

DEFINITION 1 (Max-likelihood fit). In the *max-likelihood fit* problem, we are given an arbitrary sample $S \subseteq \Re^n$ and a number $k$; we desire the Gaussian mixture with $k$ components that maximizes the likelihood of $S$.

2.2. *Classification problem.* Now we define the subcase of the learning problem when the data is assumed to arise from an unknown mixture of $k$ Gaussians, where $k$ is known.

DEFINITION 2 (Classification problem). In the *classification problem*, we are given an integer $k$, a real number $\delta > 0$ and a sample $S$ generated from an unknown mixture $F_1, F_2, \ldots, F_k$ of $k$ Gaussians in $\Re^n$, where the mixing weights $w_1, w_2, \ldots, w_k$ are also unknown. The goal is to find the "correct" labeling for each point in $S$ (up to permutation), namely to partition $S$ into $k$ subsets such that with probability at least $1 - \delta$, the partition of $S$ is exactly into the subsets of samples drawn according to each $F_i$.

Viewing the unknown mixture as a "source" we may view this as the "source learning" problem. Note that once we know the partition, we can immediately obtain estimates of the unknown Gaussians and their mixing weights.

So, the classification problem has a stronger hypothesis than the maximum-likelihood problem in that it assumes that the data came from a mixture. It also then requires the result to satisfy the stronger requirement that it is exactly the partition into the actual $S_1, S_2, \ldots, S_k$, where $S_i$ was generated according to the Gaussian $F_i$. (We abuse notation here slightly; we can only know the real $S_i$ up to a permutation of their indices. However, to avoid extra notation, we will say the partition has to be $S_1, S_2, \ldots, S_k$.)

2.3. *Our results.* Our main result concerns the classification problem. Clearly, the problem has no unique solution if the Gaussians in the mixture are allowed to be arbitrarily close to each other. We will assume a certain separation between the centers of the Gaussians. The required separation is an important consideration and will be motivated in detail in Section 3.2. Here we will just state it and mention two important features of it.

NOTATION. First, we introduce some notation which we will use throughout. We let $p_1, p_2, \ldots, p_k$ denote the (unknown) centers, respectively, of the $k$ Gaussians $F_1, F_2, \ldots, F_k$ comprising the mixture; the maximum variance of $F_i$ in any direction will be denoted $\sigma_{i,\max}$. We denote by $R_i$ the "median radius" of $F_i$; $R_i$ has the property that the $F_i$-measure of a ball of radius $R_i$ around $p_i$ is exactly $1/2$. Henceforth, we will reserve the word "radius" of a Gaussian to mean its median-radius.

Here is our formal definition of separation.

DEFINITION 3. For any $t > 0$, we say that the mixture is *t-separated* if

$$(2) \qquad \begin{aligned} |p_i - p_j|^2 &\geq -|R_i^2 - R_j^2| + 500t(R_i + R_j)(\sigma_{i,\max} + \sigma_{j,\max}) \\ &\quad + 100t^2(\sigma_{i,\max}^2 + \sigma_{j,\max}^2) \qquad \forall i, j. \end{aligned}$$

We point out here quickly two features of this definition. First, if two Gaussians $F_i, F_j$ are both spherical of the same radii ($R_i = R_j$), then the required separation is $O^*(R_i/n^{1/4})$. Second, if $F_i, F_j$ are still spherical, but if $R_j > R_i$, the term $-|R_i^2 - R_j^2|$ is negative and makes the separation required less. Indeed if $R_j = (1 + \Omega^*(1/\sqrt{n}))R_i$, then the two Gaussians $F_i, F_j$ are allowed be to concentric! The superscript $*$ on $O, \Omega$ indicates that we have omitted factors logarithmic in $n$.

THEOREM 1. *There is a polynomial-time algorithm for the classification problem. The algorithm needs to know a lower bound $w_{\min}$ on the mixing weights, and the number $s$ of sample points required is $O(n^2k^2 \log(kn^2)/(\delta^2 w_{\min}^6))$. The Gaussians may have arbitrary shape but have to be $t$-separated, where $t = O(\frac{\log s}{\delta})$.*

We also present an approximation algorithm for a special case of the max-likelihood fit problem.

THEOREM 2.    *There is a polynomial-time approximation algorithm for the max-likelihood fit problem in $\Re^n$ when the Gaussians to be fitted to the data have to be spherical of equal radii (the radius and the centers of the $k$ Gaussians have to be determined by the algorithm). There is a fixed constant $c$ such that the algorithm produces a solution whose log-likelihood is at least the best possible minus $c$.*

The algorithm of Theorem 2 is combinatorial and appears in Section 4. We note even this subcase of the maximum-likelihood fit problem is at least as hard as the clustering problem $k$-median (sum-of-squares version with Steiner nodes), which is NP-hard [8]. Indeed, our algorithm is obtained by reducing to the $k$-median algorithm of [4] (recent more efficient $k$-median algorithms would also work). We feel that this way of viewing the learning problem as an approximation problem may be fruitful in other contexts.

2.4.  *Comparison with other results.*  The algorithm in [5] makes the following assumptions: (i) all the Gaussians have the same variance–covariance matrix $\Sigma$; (ii) the maximum and minimum eigenvalues $\sigma_{\max}^2$ and $\sigma_{\min}^2$, respectively, of $\Sigma$ satisfy $\frac{\sigma_{\max}}{\sigma_{\min}} \in O(\sqrt{n}/\log k)$; (iii) the centers of any two of the $k$ Gaussians are at distance (at least) $\Omega(\sqrt{n}\sigma_{\max})$ apart.

Dasgupta and Schulman [6] showed that the EM algorithm learns (and indeed does so in just two rounds) a mixture of *spherical* Gaussians $F_1, F_2, \ldots, F_k$, where $F_i$ has radius $R_i$ (the $R_i$ may be different). They require now only a separation between centers of $F_i, F_j$ of $\Omega((R_i + R_j)/n^{1/4})$. (This amount of separation ensures among other things that the densities are "nonoverlapping"; i.e., there are $k$ disjoint balls, each containing the samples picked according to one $F_i$.)

As mentioned, our result is stronger in two ways. First, we allow Gaussians of arbitrary (and different) variance–covariance matrices and, second, we allow densities to overlap, or even be concentric. More specifically, the term $-|R_i^2 - R_j^2|$ (which is nonpositive) can make the minimum required separation negative (and so a vacuous requirement) in some cases; it allows the centers to be close (or even coincide) if the radii are very different. This allows a "large feature" to have an identifiable smaller "feature" buried inside. For the case dealt with by [6], their requirement is the same as ours (since in this case $R_i \approx \sqrt{n}\sigma_{i,\max}$) but for this term and logarithmic factors and thus their result essentially follows as a special case of ours. For the case dealt with by [5], our requirement is again weaker than that paper's but for logarithmic factors [since $\sqrt{n}\sigma_{i,\max} \in O(R_i)$]. After the first appearance of our paper [2], Vempala and Wong [25] improved the separation requirement to essentially the optimal one for the special case of spherical Gaussians: $|p_i - p_j| = \Omega((R_i + R_j)/\sqrt{n})$. Their spectral technique is entirely different from ours.

**3. Algorithm for classification problem.** First we fix notation for the rest of the section. We are given a set $S$ of samples, picked according to an unknown mixture $w_1 F_1 + w_2 F_2 + \cdots + w_k F_k$ of Gaussians $F_1, F_2, \ldots, F_k$. The known quantities are $k$ and a number $w_{\min}$ that is a lower bound on the $w_i$'s. We have to decompose $S$ as $S = S_1 \cup S_2 \cup \cdots \cup S_k$, where $S_i$ are the samples that were generated using $F_i$.

Section 3.1 describes the algorithm at an intuitive level. This description highlights the need for a "well-separated" condition on the Gaussians, which we explain in Section 3.2. The description also highlights the need for "distance concentration" results for Gaussians, which are then proved in Section 3.3. In Section 3.5 we formally describe the algorithm and prove its correctness.

3.1. *Algorithm overview.* The algorithm uses distance-based clustering, meaning that we repeatedly identify some sample point $x$ and some distance $r$, and all sample points in $B(x, r)$ all put into the same cluster. Such distance-based clustering is not new and it appears in many heuristics, including [5, 6]. The choice of $x, r$ is of course the crucial new element we provide. Since distance-based methods seem restrictive at first glance, the surprising part here is that we get provable results which subsume previously known provable results for any algorithm. This power arises from a "bootstrapping" idea, whereby we learn a little about the Gaussians from a coarse examination of the data and then bootstrap from that information to find a better clustering.

In general, distance-based clustering is most difficult when the Gaussians have different shapes and sizes, and overlap with each other (all of which we allow). It is easy to see why: a sample point from Gaussian $F_i$ might be closer to some sample points of another Gaussian $F_j$ than to all the sample points of $F_i$. One crucial insight in our algorithm is that this is unlikely to happen if we look at the Gaussian with the smallest radius in the mixture; hence we should use clustering to first identify this Gaussian, and then iterate to find the remaining Gaussians.

Now we give an overview of the algorithm. Let $F_1$ be the Gaussian of smallest radius. Using our distance-concentration results, we can argue that for any $x \in S_1$, there is an $r$ such that (i) $B(x, r) \cap S = S_1$; (ii) there is a "sizable" gap after $r$, namely, the annulus $B(x, r') \setminus B(x, r)$ for some $r'$ noticeably larger than $r$ contains no samples from any $S_j$ for $j > 1$; (iii) there is no spurious large gaps before $r$, which would confuse the algorithm.

Even after proving the above statements, the design of the algorithm is still unclear. The problem is to figure out the size of the gap between where $S_1$ ends and $S \setminus S_1$ begins, so we know when to stop. (Note: there will be gaps before $r$; the point is that they will be smaller than the one after $r$.) Our separation condition ensures that the gap between $S_1$ and the other $S_j$ is $\Omega(\sigma_{1,\max})$; so we need an estimate of $\sigma_{1,\max}$. We get such an estimate by bootstrapping. We show that if we have any fraction $f$ of the samples in $S_1$, then we may estimate $\sigma_{1,\max}$ to a factor of $O(1/f^2)$ with high probability. We use this to get a rough estimate $\beta$ of $\sigma_{1,\max}$. Using $\beta$, we increment the radius in steps which are guaranteed to be less than $\sigma_{1,\max}$ (which ensures that we do not step over the "gap" into another $S_j$) until we observe a gap; by then, we have provably picked up *most* of $S_1$. Now we use this to better estimate $\sigma_{1,\max}$ and then incrementing the radius by another $\Omega(\sigma_{1,\max})$, we capture all of $S_1$. (The guaranteed gap ensures that we do not get any points from any other Gaussian while we increment the radius.)

To make all the above ideas rigorous, we need appropriate distance-concentration results which assert that the distance between certain pairs of sample points considered as a random variable is concentrated around a certain value. Some distance-concentration results—at least for spherical or spherelike Gaussians—were known prior to our work, showing a sharp concentration around the mean or median. However, for the current algorithm we also need concentration around values that are quite far from the mean or median. For example, to show the nonexistence of "spurious gaps," we have to show that if a ball of radius $r$ centered at a sample point $x \in S_i$ has $F_i$-measure, say, exactly $1/4$, then, for a small $\delta > 0$, the ball of radius $r + \delta$ with $x$ as center has $F_i$-measure at least $0.26$. If such a result failed to hold, then we might see a "gap" (an annulus with no sample points) and falsely conclude that we have seen all of $S_i$. Such concentration results (around values other than the median or mean) are not in general provable by the traditional moment-generating function approach. We introduce a new approach in this context: isoperimetric inequalities (see Theorem 3). Our method does not always prove the tightest possible concentration results, but is more general. For example, one may derive weaker concentration results for general log-concave densities via this method (see [17]).

3.2. *Separation condition and its motivation.* Now we motivate our separation condition, which is motivated by the exigencies of distance-based clustering. Consider the very special case of spherical Gaussians $F_i, F_j$ with $R_i = R_j$. Suppose $x, x', y$ are independent samples, $x, x'$ picked according to $F_i$ and $y$ according to $F_j$. Lemma 5 will argue, that with high probability [we will use $\approx$ here to mean that the two sides differ by at most $O(R_i^2/\sqrt{n})$],

$$|x - p_i|^2 \approx R_i^2$$

and similar concentration results for $|x' - p_i|, |y - p_j|$. It is an intuitive fact that $x - p_i, x' - p_i, p_i - p_j, y - p_j$ will all be pairwise nearly orthogonal (a sample from a spherical Gaussian is almost orthogonal to any fixed direction with high probability). So, one can show that

$$(3) \qquad |x - x'|^2 \approx |x - p_i|^2 + |p_i - x'|^2 \approx 2R_i^2,$$

$$(4) \qquad |x - y|^2 \approx |x - p_i|^2 + |p_i - p_j|^2 + |y - p_j|^2 \approx 2R_i^2 + |p_i - p_j|^2.$$

(The first assertion is proved rigorously in greater generality in Lemma 7 and the second one in Lemma 8.) Thus, it is clear that if $|p_i - p_j|^2$ is at least $\Omega(R_i^2/\sqrt{n})$, then with high probability $|x - x'|$ and $|x - y|$ will lie in different ranges. (Aside: One can also show a sort of converse with different constants, since the concentration results we get are qualitatively tight. However, we will not establish this, since it is not needed.) This intercenter separation then is

$$O(R_i/n^{1/4}).$$

Our separation condition for this case is indeed this quantity, up to a factor $\log n$. A weaker separation condition would be to require a separation of $\Omega(R_i/\sqrt{n})$; at this separation, one can still show that with high probability the hyperplane orthogonal to the line joining the centers at its midpoint has all the samples of one Gaussian on one side and the samples of the other Gaussian on the other side. An algorithm to learn under this condition would be a stronger result than our distance-based algorithm in this case. Since the conference version of our paper appeared, Vempala and Wang [25] have indeed developed a learning algorithm under this weaker separation for the case of spherical Gaussians using spectral techniques.

3.3. *Concentration results using isoperimetric inequalities.* Suppose we have some probability density $F$ in $\Re^n$ and a point $x$ in space. For proving distance concentration results, we would like to measure the rate of growth or decline of $F(B(x, r))$ as a function of $r$. This will be provided by the isoperimetric inequality (see Corollary 4).

THEOREM 3 [14]. *Let $F(x) = e^{-x^T A^{-1} x} g(x)$ be a function defined on $\Re^n$, where $A$ is a positive definite matrix whose largest eigenvalue is $\sigma_{\max}^2$ and $g(x)$ is any positive real-valued log-concave function. Suppose $\nu$ is a positive real and we have a partition of $\Re^n$ into three sets $K_1, K_2, K_3$ so that, for all $x \in K_1, y \in K_2$, we have $|x - y| \geq \nu$. Then*

$$\int_{K_3} F(x)\, dx \geq \frac{2\nu e^{-\nu^2}}{\sqrt{\pi}} \frac{1}{\sigma_{max}} \min\left( \int_{K_1} F(x)\, dx, \int_{K_2} F(x)\, dx \right).$$

The phrase "isoperimetric inequality" has come to mean a lower bound on the surface area of a set in terms of its volume. If $K_1$ is fixed and we define $K_3$ to be the set of points not in $K_1$ which are at distance at most $\nu$ from some point in $K_1$ and define $K_2 = \Re^n \setminus (K_1 \cup K_3)$, then as $\nu$ goes to zero, $K_3$ tends to the boundary surface of $K_1$ and the above theorem can be shown to yield a lower bound on the surface integral of $F$ over this surface. We will make this connection rigorous below for the context we need. Such isoperimetric inequalities for general log-concave measures over multidimensional sets were first proved for use in establishing rapid convergence to the steady state of certain Markov chains for approximating volumes of convex sets and for sampling according to log-concave measures [1, 9, 16]. The proof of Theorem 3 uses a specialization of the above techniques to the case of Gaussians, where we get better results.

COROLLARY 4.   *We borrow notation from Theorem* 3 *and also assume that* $F(\Re^n) = 1$:

(i) *If a ball* $B(x, r)$ *has* $F(B(x, r)) \leq 1/2$, *then*
$$\frac{d(\ln(F(B(x, r))))}{dr} \geq \frac{2}{\sqrt{\pi} \sigma_{\max}}.$$

(ii) *If a ball* $B(x, r)$ *has* $F(B(x, r)) \geq 1/2$, *then*
$$\frac{d(\ln(1 - F(B(x, r))))}{dr} \leq \frac{-2}{\sqrt{\pi} \sigma_{\max}}.$$

REMARK.   The corollary says that $\ln(F(B(x, r)))$ grows at a rate of $\Omega(1/\sigma_{\max})$ until $F(B(x, r))$ is $1/2$, and then $\ln(1 - F(B(x, r)))$ declines at a rate of $\Omega(1/\sigma_{\max})$. Intuitively, it is easy to see that this would lead to distance concentration results since once we increase (decrease) $r$ by $O(\sigma_{\max})$ from its median value, the mass outside $B(x, r)$ [inside $B(x, r)$] is small. The first lemma below (Lemma 5) is derived exactly on these lines; the subsequent three Lemmas 6–8 discuss the distances between different samples from the same and from different Gaussians.

PROOF OF COROLLARY 4.   Let $\nu$ be an infinitesimal. Then
$$\frac{d(\ln(F(B(x, r))))}{dr} = \frac{1}{F(B(x, r))} \frac{d(F(B(x, r)))}{dr}$$
$$= \lim_{\nu \to 0} \frac{1}{\nu F(B(x, r))} [F(B(x, r + \nu)) - F(B(x, r))].$$

Now we let $K_1 = B(x, r)$ and $K_2 = \Re^n \setminus B(x, r + \nu)$ and apply the theorem above to get the first assertion of the corollary. The second assertion follows similarly.

$\square$

LEMMA 5. *Suppose $F$ is a general Gaussian in $\Re^n$ with maximum variance in any direction $\sigma$, radius $R$ and center $p$. Then, for any $t > 0$, we have*

$$F(\{x : R - t\sigma \leq |x - p| \leq R + t\sigma\}) \geq 1 - e^{-t}.$$

PROOF. For any $\gamma > 0$, let $F(B(p, \gamma)) = g(\gamma)$. Then, for $\gamma < R$, we have by Corollary 4 that

$$\frac{d \ln(g(\gamma))}{d\gamma} \geq \frac{1}{\sigma}.$$

Integrating from $\gamma$ to $R$, we get that

$$F(B(p, \gamma)) \leq \tfrac{1}{2} e^{-(R-\gamma)/\sigma}.$$

For $\gamma > R$, isoperimetry implies that

$$\frac{d(\ln(1 - g(\gamma)))}{d\gamma} \leq \frac{-1}{\sigma}.$$

Again integrating from $R$ to $\gamma$, we get $1 - g(\gamma) \leq (1/2)e^{-(\gamma - R)/\sigma}$. Combining the two, the lemma follows. □

LEMMA 6. *Let $F, p, R, \sigma$ be as in Lemma 5 and suppose $z$ is any point in space. Let $t \geq 1$. If $x$ is picked according to $F$, we have that, with probability at least $1 - 2e^{-t}$,*

(5)
$$\begin{aligned}
(R + t\sigma)^2 &+ |z - p|^2 + 2\sqrt{2}\sqrt{t}|z - p|\sigma \\
&\geq |x - z|^2 \\
&\geq ((R - t\sigma)^+)^2 + |z - p|^2 - 2\sqrt{2}\sqrt{t}|z - p|\sigma,
\end{aligned}$$

*where $(R - t\sigma)^+$ is $R - t\sigma$ if this quantity is positive and 0 otherwise.*

PROOF. We have

(6)
$$\begin{aligned}
|x - z|^2 &= ((x - p) + (p - z)) \cdot ((x - p) + (p - z)) \\
&= |x - p|^2 + |p - z|^2 + 2(x - p) \cdot (p - z).
\end{aligned}$$

Now $2(x - p) \cdot (p - z)$ is a normal random variable with mean 0 and variance at most $4|p - z|^2\sigma^2$, so the probability that $|2(x - p) \cdot (p - z)|$ is greater than $2\sqrt{2}\sqrt{t}|z - p|\sigma$ is at most $e^{-t}$. From Lemma 5, we have that $R - t\sigma \leq |x - p| \leq R + t\sigma$ with probability at least $1 - e^{-t}$. Combining these two facts, the current lemma follows. □

LEMMA 7. *Suppose $F, p, R, \sigma$ as in Lemma 5. Suppose $x, y$ are independent samples each picked according to $F$. Then for any $t \geq 1$, with probability at least $1 - 3e^{-t}$, we have*

$$2R^2 - 8t\sigma R \leq |x - y|^2 \leq 2(R + 2t\sigma)^2.$$

PROOF. We may assume that $x$ is picked first and then $y$ (independently). Then from Lemma 5, with probability $1 - e^{-t}$, we have $R - t\sigma \leq |x - p| \leq R + t\sigma$. From Lemma 6 (once $x$ is already picked), with probability at least $1 - 2e^{-t}$, we have $(R + t\sigma)^2 + |x - p|4\sqrt{t}\sigma + |x - p|^2 \geq |x - y|^2 \geq R^2 - 2Rt\sigma - 4|x - p|\sqrt{t}\sigma + |x - p|^2$. Both conclusions hold with probability at least $1 - 3e^{-t}$, whence we get

$$|x - y|^2 \leq (R + t\sigma)^2 + 4t\sigma(R + t\sigma) + (R + t\sigma)^2 \leq 2(R + 2t\sigma)^2.$$

For the lower bound on $|x-y|^2$, first note that if $R \leq 4t\sigma$, then $2R^2 - 8t\sigma \leq 0$, so the lower bound is obviously valid. So we may assume that $R > 4t\sigma$. Thus, $|x - p| \geq 3t\sigma$ and under this constraint, $|x - p|^2 - 4\sqrt{t}\sigma|x - p|$ is an increasing function of $|x - p|$. So, we get

$$|x - y|^2 \geq R^2 - 2Rt\sigma - 4t\sigma(R - t\sigma) + (R - t\sigma)^2,$$

which yields the lower bound claimed.   □

LEMMA 8.   *Let $t \geq 1$. If $x$ is a random sample picked according to $F_i$ and $y$ is picked independently according to $F_j$, with $F_i, F_j$ satisfying the separation condition (2), then, with probability at least $1 - 6e^{-t}$, we have*

$$\begin{align}
(7) \qquad |x - y|^2 &\geq 2\min(R_i^2, R_j^2) + 60t(\sigma_{i,\max} + \sigma_{j,\max})(R_i + R_j) \\
&\quad + 30t^2(\sigma_{i,\max}^2 + \sigma_{j,\max}^2).
\end{align}$$

PROOF.   Assume without loss of generality that $R_i \leq R_j$. Applying Lemma 6, we get that, with probability at least $1 - 2e^{-t}$, we have

$$(8) \quad |y - p_i|^2 \geq R_j^2 - 2t\sigma_{j,\max}R_j + |p_i - p_j|^2 - 2\sqrt{2}\sqrt{t}|p_i - p_j|\sigma_{j,\max}.$$

CLAIM 1.

$$(9) \quad |y - p_i|^2 \geq R_i^2 + 154t(\sigma_{i,\max} + \sigma_{j,\max})(R_i + R_j) + 30t^2(\sigma_{i,\max}^2 + \sigma_{j,\max}^2).$$

PROOF.   *Case* 1. $R_j^2 \geq R_i^2 + 250t(R_i + R_j)(\sigma_{i,\max} + \sigma_{j,\max}) + 30t^2(\sigma_{i,\max}^2 + \sigma_{j,\max}^2)$. Note that $|p_i - p_j|^2 - 4\sqrt{t}\sigma_{j,\max}|p_i - p_j| + 4t\sigma_{j,\max}^2 \geq 0$. So, $|p_i - p_j|^2 - 4\sqrt{t}\sigma_{j,\max}|p_i - p_j| \geq -4t\sigma_{j,\max}^2$. Plugging this into (8), and using the case assumption, we get

$$\begin{align}
|y - p_i|^2 &\geq R_i^2 + 250t(R_i + R_j)(\sigma_{j,\max} + \sigma_{i,\max}) \\
&\quad + 30t^2(\sigma_{j,\max}^2 + \sigma_{i,\max}^2) - 2t\sigma_{j,\max}R_j - 4t\sigma_{j,\max}^2.
\end{align}$$

It is easy to see that $R_j \geq (2/3)\sigma_{j,\max}$—this is because $R_j^2$ is clearly at least the median value of $(u \cdot x)^2$ under $F_j$, where $u$ is the direction achieving $\sigma_{j,\max}$; now it is easy see that, for the one-dimensional Gaussian $u \cdot x$, the median value of $(u \cdot x)^2$ is at least $2/3$ times the variance by direct calculation. Plugging $\sigma_{j,\max}^2 \le (3/2)R_j\sigma_{j,\max}$ into the above inequality, we easily get the claim.

*Case* 2. $R_j^2 < R_i^2 + 250t(R_i + R_j)(\sigma_{i,\max} + \sigma_{j,\max}) + 30t^2(\sigma_{i,\max}^2 + \sigma_{j,\max}^2)$. Then by the separation condition, we have

$$|p_i - p_j|^2 \ge 250t(R_i + R_j)(\sigma_{i,\max} + \sigma_{j,\max}) + 70t^2(\sigma_{i,\max}^2 + \sigma_{j,\max}^2).$$

Now, since $|p_i - p_j|^2 - 2\sqrt{2}\sqrt{t}\sigma_{j,\max}|p_i - p_j|$ is an increasing function of $|p_i - p_j|$ for $|p_i - p_j| \ge 2\sqrt{2}\sqrt{t}\sigma_{j,\max}$, we have

$$|p_i - p_j|^2 - 2\sqrt{2}\sqrt{t}\sigma_{j,\max}|p_i - p_j|$$
$$\ge 250t(R_i + R_j)(\sigma_{i,\max} + \sigma_{j,\max}) + 70t^2(\sigma_{i,\max}^2 + \sigma_{j,\max}^2)$$
$$\quad - 2\sqrt{2}\sqrt{t}\sigma_{j,\max}(16\sqrt{t}\sqrt{R_j + R_i}\sqrt{\sigma_{j,\max} + \sigma_{i,\max}} + 9t(\sigma_{j,\max} + \sigma_{i,\max}))$$
$$\ge 156t(R_i + R_j)(\sigma_{i,\max} + \sigma_{j,\max}) + 34t^2(\sigma_{i,\max}^2 + \sigma_{j,\max}^2),$$

using the inequality $\sqrt{a+b} \le \sqrt{a} + \sqrt{b}$, $\forall a, b \ge 0$ and observing that $\sigma_{j,\max} \le (3/2)R_j$ and $\sigma_{j,\max}(\sigma_{j,\max} + \sigma_{i,\max}) \le \sqrt{2}(\sigma_{i,\max}^2 + \sigma_{j,\max}^2)$.

Putting this into (8), we get

$$|y - p_i|^2 \ge R_i^2 - 2t\sigma_{j,\max}R_j + 156t(R_i + R_j)(\sigma_{i,\max} + \sigma_{j,\max})$$
$$\quad + 34t^2(\sigma_{i,\max}^2 + \sigma_{j,\max}^2),$$

which yields the claim in this case. $\square$

Imagine now $y$ already having been picked and $x$ being picked independently of $y$. Applying Lemma 6, we get that, with probability at least $1 - 2e^{-t}$, we have (again using the inequality, $\sqrt{a+b} \le \sqrt{a} + \sqrt{b}$, $\forall a, b \ge 0$)

$$|x - y|^2 \ge R_i^2 - 2R_i t\sigma_{i,\max} + |y - p_i|^2 - 2\sqrt{2}\sqrt{t}\sigma_{i,\max}|y - p_i|$$
$$\ge R_i^2 - 2R_i t\sigma_{i,\max} + R_i^2 + 154t(\sigma_{i,\max} + \sigma_{j,\max})(R_i + R_j)$$
$$\quad + 30t^2(\sigma_{i,\max}^2 + \sigma_{j,\max}^2)$$
$$\quad - 2\sqrt{2}\sqrt{t}\sigma_{i,\max}(R_i + 13\sqrt{t}\sqrt{\sigma_{i,\max} + \sigma_{j,\max}}\sqrt{R_i + R_j}$$
$$\quad + \sqrt{30}t(\sigma_{i,\max} + \sigma_{j,\max}))$$

[because under condition (9), $|y - p_i|^2 - 2\sqrt{2}\sqrt{t}\sigma_{i,\max}|y - p_i|$ is an increasing function of $|y - p_i|$]

$$|x - y|^2 \ge 2R_i^2 + (154 - 2 - 2\sqrt{2}(1 + 13\sqrt{3/2} + 1.5\sqrt{30}))$$
$$\quad \times t(\sigma_{i,\max} + \sigma_{j,\max})(R_i + R_j) + 30t^2(\sigma_{i,\max}^2 + \sigma_{j,\max}^2),$$

from which the lemma follows. $\square$

3.4. *Warm-up*: *case of spherical Gaussians.* As a consequence of our concentration results we first present our algorithm for the simple case when all the $F_i$ are spherical. In this case, $\sigma_{i,\max} \approx R_i/\sqrt{n}$, where the error is small enough that our calculations below are valid. Choosing $t = \Omega(\log(|S|/\delta))$ as before, it is easy to see from the distance concentration results that, with high probability,

$$(10) \quad |x - y|^2 \in \left[ 2R_i^2\left(1 - \frac{4t}{\sqrt{n}}\right), 2R_i^2\left(1 + \frac{5t}{\sqrt{n}}\right) \right] \qquad \forall\, x, y \in S_i,\, \forall\, i,$$

and by appropriately choosing the constant in the definition of $t$-separation we can also ensure that with high probability there is a positive constant $c' > 12$ such that

$$(11) \quad |x - y|^2 \geq 2\min(R_i^2, R_j^2) + \frac{c't(R_i + R_j)^2}{\sqrt{n}} \qquad \forall\, x \in S_i,\, \forall\, y \in S_j,\, \forall\, i \neq j.$$

For each pair $x, y \in S$ find $|x - y|$ and suppose $x_0, y_0$ is a pair (there may be several) at the minimum distance. Then from (10) and (11) it follows that if $x_0 \in S_i$, then for all $y \in S_i$, $|x - y| \leq (1 + \frac{3t}{\sqrt{n}})|x_0 - y_0|$ and furthermore, for all $z \in S \setminus S_i$, $|x - z| \geq (1 + \frac{6t}{\sqrt{n}})|x_0 - y_0|$. So, we may identify $S_i$ by $S \cap B(x_0, |x_0 - y_0|(1 + \frac{3t}{\sqrt{n}}))$. Having thus found $S_i$, we may peel it off from $S$ and repeat the argument. The important thing here is that we can estimate the radius of the ball—namely, $|x_0 - y_0|(1 + \frac{3t}{\sqrt{n}})$—from observed quantities; this will not be so easily the case for general Gaussians.

3.5. *The general case.* Now we consider the case when the Gaussians may not be spherical or even spherelike. Let $\delta > 0$ be the probability of failure allowed. We are given a set of samples $S$ drawn according to an unknown mixture of Gaussians $w_1 F_1 + w_2 F_2 + \cdots + w_k F_k$; but we are given a $w_{\min} > 0$ with $w_i \geq w_{\min}$ for all $i$. We assume that $|S| \geq 10^7 n^2 k^2 \log(kn^2)/(\delta^2 w_{\min}^6)$. In what follows, we choose

$$t = \frac{100 \log |S|}{\delta}.$$

THE ALGORITHM.    Initialization: $T \leftarrow S$.

1. Let $\alpha > 0$ be the smallest value such that a ball $B(x, \alpha)$ of radius $\alpha$ centered at some point in $x \in T$ has at least $3w_{\min}|S|/4$ points from $T$. (This will identify a Gaussian $F_i$ with approximately the least radius.)
2. Find the maximum variance of the set $Q = B(x, \alpha) \cap T$ in any direction. That is, find

$$\beta = \max_{w\,:\,|w|=1} \frac{1}{|Q|} \sum_{y \in Q} \left( w \cdot y - w \cdot \left( \frac{1}{|Q|} \sum_{z \in Q} z \right) \right)^2.$$

(This $\beta$ is our first estimate of $\sigma_{\max}$. Note that computing $\beta$ is an eigenvalue computation, and an approximate eigenvalue suffices.)

3. Let $\nu = \sqrt{\frac{w_{\min}\beta}{8}}$. (We will later show that $\nu \leq \sigma_{\max}$; so increasing the radius in steps of $\nu$ ensures that we do not miss the "gap" between the $S_i$ that $x$ belongs to and the others.) Find the least positive integer $s$ such that (we will later prove that such a $s$ exists)

$$B(x, \alpha + s\nu) \cap T = B(x, \alpha + (s-1)\nu) \cap T.$$

4. Let $Q' = B(x, \alpha + s\nu) \cap T$. As in step 3, find the maximum variance $\beta'$ of the set $Q'$ in any direction. (We will prove that this $\beta'$ gives a better estimate of $\sigma_{\max}$.)

5. Remove $B(x, \alpha + s\nu + 3\sqrt{\beta'}(\log|S| - \log\delta + 1)) \cap T$ from $T$. (We will show that the set so removed is precisely one of the $S_i$.)

6. Repeat until $T$ is empty.

REMARK 1. Ball $B(x, \alpha + s\nu)$ will be shown to contain all but $w_{\min}/(10w_i)$ of the mass of the Gaussian $F_i$ we are dealing with; the bigger radius of $B(x, \alpha + s\nu + 3\sqrt{\beta'}(\log|S| - \log\delta + 1))$ will be shown to include all but $\delta/(10|S|^2)$ of the mass of $F_i$. This will follow using isoperimetry. Then we may argue that with high probability all of $S_i$ is now inside this ball. An easier argument shows that none of the other $S_j$ intersect this ball.

Now we prove why this works as claimed. Let $\delta > 0$ be the probability of failure allowed. Recall that

$$t = \frac{100\log|S|}{\delta}.$$

We will now show using the distance-concentration results that several desirable events [described below in (12)–(18)] happen, each with probability at least $1 - \frac{\delta}{10}$. *We will assume from now on that conditions* (12)–(18) *hold after allowing for the failure probability of at most* $7\delta/10$. The bottom line is that the sample is very likely to represent the mixture accurately: the component Gaussians are represented essentially in proportion to their mixing weights; the number of samples in every sphere and half-space is about right and so forth.

First, since $|S_i|$ can be viewed as the sum of $|S|$ Bernoulli independent 0–1 random variables, where each is 1 with probability $w_i$, we have [using standard results, e.g., Hoeffding's inequality, which asserts that for $s$ i.i.d. Bernoulli random variables $X_1, X_2, \ldots, X_s$ with $\text{Prob}(X_i = 1) = q$, for all real numbers $\alpha > 0$, $\text{Prob}(|\sum_{i=1}^{s} X_i - sq| \geq \alpha) \leq 2e^{-\alpha^2 q/4s}$] that, with probability at least $1 - \delta/10$,

(12)                    $1.1w_i|S| \geq |S_i| \geq 0.9w_i|S| \qquad \forall i.$

For each $i, 1 \le i \le k$, and each $x \in S_i$, let $\eta(x)$ be the least positive real number such that

$$F_i(B(x, \eta(x))) \ge 1 - \frac{\delta}{10|S|^2}.$$

Now, we assert that, with probability at least $1 - \frac{\delta}{10}$,

$$(13) \qquad \forall\, i, 1 \le i \le k,\, \forall\, x \in S_i \qquad S_i \subseteq B(x, \eta(x)).$$

To see this, focus attention first on one particular $x \in S$, say $x \in S_i$. We may imagine picking $x$ as the first sample in $S$ and then independently picking the others. Then since $x$ is fixed, $\eta(x)$ and $B(x, \eta(x))$ are fixed; so from the lower bound on $F_i(B(x, \eta(x)))$, it follows that $\mathrm{Prob}(S_i \subseteq B(x, \eta(x))) \ge 1 - \delta/(10|S|)$. From this we get (13).

We have from Lemma 7 that, with probability at least $1 - \delta/10$, the following is true for each $i, 1 \le i \le k,\, \forall\, x, y \in S_i$:

$$(14) \qquad 2R_i^2 - 8t\sigma_{i,\max}R_i \le |x - y|^2 \le 2(R_i + 2t\sigma_{i,\max})^2.$$

Further, from Lemma 8, we have that, with probability at least $1 - \frac{\delta}{10}$,

$$(15) \qquad \begin{aligned} &\forall\, i, j \le k, i \ne j,\, \forall\, x \in S_i,\, \forall\, y \in S_j \\ &|x - y|^2 \ge 2\min(R_i^2, R_j^2) + 60t(R_i + R_j)(\sigma_{i,\max} + \sigma_{j,\max}) \\ &\qquad + 30t^2(\sigma_{i,\max}^2 + \sigma_{j,\max}^2). \end{aligned}$$

Next, we wish to assert that certain spherical annuli centered at sample points have roughly the right number of points. Namely,

$$(16) \qquad \begin{aligned} &\forall\, i,\, \forall\, x, y, z \in S_i \text{ letting } A = B(x, |x - y|) \setminus B(x, |x - z|) \\ &\text{we have } \left| \frac{|S_i \cap A|}{|S_i|} - F_i(A) \right| \le \frac{w_{\min}^{5/2}}{160}. \end{aligned}$$

We will only sketch the routine argument for this. First, for a particular triple $x, y, z$ in some $S_i$, we may assume that these points $x, y, z$ were picked first and the other points of the sample are then being picked independently. So for the other points, the annulus is a fixed region in space. Then we may view the rest as Bernoulli trials and apply Hoeffding's inequality. The above follows from the fact that the Hoeffding upper bound multiplied by $|S|^3$ (the number of triples) is at most $\delta/10$.

Next, we wish to assert that every half-space in space contains about the correct number of sample points. For this, we use a standard Vapnik–Chervonenkis (VC) dimension argument [24]. They define a fundamental notion called VC dimension (which we do not define here). If a (possibly infinite) collection $\mathcal{C}$ of subsets of $\Re^n$ has VC dimension $d$ and $\mathcal{D}$ is an arbitrary probability distribution on $\Re^n$, then for any $\rho, \varepsilon > 0$ and for a set of

$$\frac{4}{\varepsilon} \log \frac{2}{\rho} + \frac{8d}{\varepsilon} \log \frac{8d}{\varepsilon}$$

independent identically distributed samples drawn according to $\mathcal{D}$, we have that with probability at least $1 - \rho$, for *every* $H \in \mathcal{C}$, the fraction of samples that lie in $H$ is between $\mathcal{D}(H) - \varepsilon$ and $\mathcal{D}(H) + \varepsilon$.

In our case, $\mathcal{C}$ consists of half-spaces; it is well known that the VC dimension of half-spaces in $\Re^n$ is $n$. We consider each component $F_i$ of our mixture in turn as $\mathcal{D}$. We have drawn a sample of size $|S_i|$ from $F_i$. Applying the VC dimension argument for each $i$, with $\rho = \delta/10k$ and $\varepsilon = w_{\min}/100$, and then using the union bound, we conclude that, with probability at least $1 - \delta/10$, the sample satisfies

(17)    $\forall\, i, \forall$ half-spaces $H$
$$||S_i \cap H| - |S_i| F_i(H)| \leq w_{\min} |S_i| / 100.$$

From Lemma 12 (to come), it follows that, with probability at least $1 - \delta/10$, we have

(18)  $\forall$ unit length vectors $w$, $\forall\, i$    $\dfrac{1}{|S_i|} \displaystyle\sum_{x \in S_i} (w \cdot (x - p_i))^2 \leq 2\sigma_{i,\max}^2.$

LEMMA 9.  *Each execution of steps* 1–5 *removes precisely one of the* $S_i$.

PROOF.  The lemma will be proved by induction on the number of executions of the loop. Suppose we have finished $l - 1$ executions and are starting on the $l$th execution.

Let $P$ be the set of $j$ such that $S_j$ has not yet been removed. (By the inductive assumption at the start of the loop, $T$ is the union of $S_j$, $j \in P$.)

LEMMA 10.  *Suppose* $x \in S$ *is the center of the ball* $B(x, \alpha)$ *found in the* $l$th *execution of step* 1 *of the algorithm and suppose* $x$ *belongs to* $S_i$ ($i$ *unknown to us*). *Then*

(19)    $B(x, \alpha) \cap S \subseteq S_i,$

(20)    $|x - y|^2 \geq 2R_i^2 + 50t(\sigma_{i,\max} + \sigma_{j,\max})(R_i + R_j)$
$\qquad\qquad + 20t^2(\sigma_{i,\max}^2 + \sigma_{j,\max}^2) \qquad \forall y \in S_j, \, \forall j \neq i, j \in P.$

PROOF.  For any $j \in P$, and all $y, z \in S_j$, we have from (14) that $|z - y|^2 \leq 2(R_j + 2t\sigma_{j,\max})^2$. Thus, a ball of radius $\sqrt{2}(R_j + 2\sigma_{i,\max})$ with $y$ as center would qualify in step 1 of the algorithm by (12). So, by definition of $\alpha$ in that step, we must have

$$(21) \qquad \alpha \le \sqrt{2}(R_j + 2t\sigma_{j,\max}) \qquad \forall\, j \in P.$$

If now $B(x, \alpha)$ contains a point $z$ from some $S_j$, $j \ne i$, by the inductive assumption in Lemma 9, we must have $j \in P$. Then by (15) we have

$$\alpha^2 \ge 2\min(R_i^2, R_j^2) + 60t(R_i + R_j)(\sigma_{i,\max} + \sigma_{j,\max}) + 30t^2(\sigma_{i,\max}^2 + \sigma_{j,\max}^2),$$

which contradicts (21) [noting that (21) must hold for both $i, j$]. This proves (19).

Now, from the lower bound of (14), it follows that

$$\alpha^2 \ge 2R_i^2 - 8R_i\sigma_{i,\max}t.$$

So from (21) it follows that

$$2R_j^2 \ge 2R_i^2 - 8t(R_i + R_j)(\sigma_{i,\max} + \sigma_{j,\max}) - 8t^2\sigma_{j,\max}^2 \qquad \forall\, j \in P.$$

Thus from (15), we get that, for $y \in S_j, j \ne i$,

$$
\begin{aligned}
(22) \quad |x - y|^2 &\ge 2R_i^2 - 8t(R_i + R_j)(\sigma_{i,\max} + \sigma_{j,\max}) - 8t^2\sigma_{j,\max}^2 \\
&\quad + 60t(R_i + R_j)(\sigma_{i,\max} + \sigma_{j,\max}) + 30t^2(\sigma_{i,\max}^2 + \sigma_{j,\max}^2),
\end{aligned}
$$

from which (20) follows.  $\square$

Now we can show that $\beta$ is a rough approximation to $\sigma_{i,\max}$.

CLAIM 2.   *The $\beta, Q$ computed in step 2 of the algorithm satisfy*

$$\frac{2|S_i|}{|Q|}\sigma_{i,\max}^2 \ge \beta \ge \frac{|Q|^2}{4|S_i|^2}\sigma_{i,\max}^2.$$

PROOF.   For any unit length vector $w$, we have, by (18),

$$\sum_{x \in Q}(w \cdot (x - p_i))^2 \le \sum_{x \in S_i}(w \cdot (x - p_i))^2 \le 2|S_i|\sigma_{i,\max}^2.$$

Since this holds for every $w$, and the second moment about the mean is less than or equal to the second moment about $p_i$, we have that $\beta \le 2\frac{|S_i|}{|Q|}\sigma_{i,\max}^2$. This proves the upper bound on $\beta$.

Let $u$ be the direction of the maximum variance of $F_i$. We wish to assert that the variance of $Q$ along $u$ is at least $|Q|\sigma_{i,\max}/|S_i|$. To this end, first note that, for any reals $\gamma_1, \gamma_2$, with $\gamma_1 > 0$, we have

$$
\begin{aligned}
\mathrm{Prob}_{F_i}(\gamma_2 - \gamma_1 \le x \cdot u \le \gamma_2 + \gamma_1) &= \frac{1}{2\sqrt{\pi}\sigma_{i,\max}}\int_{\gamma_2 - \gamma_1}^{\gamma_2 + \gamma_1} e^{-r^2/2\sigma_{i,\max}^2}\, dr \\
&\le \frac{\gamma_1}{\sqrt{\pi}\sigma_{i,\max}}.
\end{aligned}
$$

Let $\gamma_2 = \frac{1}{|Q|}\sum_{x \in Q}(u \cdot x)$ and let $\gamma_1 = \frac{|Q|}{|S_i|}\sigma_{i,\max}$. Then the strip $H = \{x : \gamma_2 - \gamma_1 \le u \cdot x \le \gamma_2 + \gamma_1\}$ satisfies $F_i(H) \le \gamma_1/(\sqrt{\pi}\sigma_{i,\max})$. So, by (17),

$$|S_i \cap H| \le |S_i|F_i(H) + \frac{w_{\min}|S_i|}{100} \le 3\frac{|Q|}{4} \qquad \text{using } |Q| \ge \tfrac{3}{4}w_{\min}|S|.$$

So, we have that

$$\frac{1}{|Q|}\sum_{x \in Q}(u \cdot x - \gamma_2)^2 \ge \frac{1}{|Q|}\frac{|Q|}{4}\frac{|Q|^2}{|S_i|^2}\sigma_{i,\max}^2 = \frac{1}{4}\sigma_{i,\max}^2\frac{|Q|^2}{|S_i|^2},$$

from which the lower bound on $\beta$ obviously follows.  □

COROLLARY 11.  *The $\beta$ computed in step 2 of the algorithm satisfies*

$$\frac{4}{w_{\min}}\sigma_{i,\max}^2 \ge \beta \ge \frac{1}{8}w_{\min}^2\sigma_{i,\max}^2.$$

PROOF.   Since $|Q| \ge 3w_{\min}|S|/4$, Claim 2 implies the corollary.   □

From (14) we have

$$\forall y \in S_i \qquad |x - y|^2 \le 2R_i^2 + 4tR_i\sigma_{i,\max} + 4t^2\sigma_{i,\max}^2.$$

From (20) we have

$$\forall z \in \bigcup_{j \in P\setminus\{i\}} S_j \qquad |x - z|^2 x \ge 2R_i^2 + 50t(\sigma_{i,\max} + \sigma_{j,\max})(R_i + R_j)$$

$$+ 20t^2(\sigma_{i,\max}^2 + \sigma_{j,\max}^2).$$

Thus, there exists an annulus of size (where the size of an annulus denotes the difference in radii between the outer and inner balls) $\sigma_{i,\max}$ around $x$ with no sample points in it. Since we are increasing the radius in steps of $\nu$ which is at most $\sigma_{i,\max}$ (Corollary 11) there is some $s$ in step 3 of the algorithm. Also, we have

$$B(x, \alpha + s\nu) \cap S \subseteq S_i.$$

The trouble of course is that such a gap may exist even inside $S_i$, so $B(x, \alpha + s\nu)$ may not contain all of $S_i$. To complete the induction we have to argue that steps 4 and 5 will succeed in identifying *every* point of $S_i$. For $\gamma > 0$, let

$$g(\gamma) = F_i(B(x, \gamma)).$$

From $B(x, \alpha + s\nu) \cap T = B(x, \alpha + (s-1)\nu) \cap T$ (see step 3 of the algorithm), we get using (16) that

$$g(\alpha + s\nu) - g(\alpha + (s-1)\nu) \le \frac{w_{\min}^{5/2}}{160}.$$

Since $\nu = \sqrt{\frac{w_{\min}\beta}{8}}$, we get using Corollary 11's lower bound on $\beta$ that there exists a $\gamma' \in [(s-1)\nu + \alpha \quad s\nu + \alpha]$ with

$$\left(\frac{dg(\gamma)}{d\gamma}\right)_{\gamma=\gamma'} \leq \frac{w_{\min}^{5/2}}{160\nu} \leq \frac{w_{\min}}{20\sigma_{i,\max}}.$$

(If not, integration would contradict the previous inequality.) Thus isoperimetry (Corollary 4) implies that

$$g(\alpha + s\nu) \geq 1 - \frac{w_{\min}}{10} \quad \text{or} \quad g(\alpha + (s-1)\nu) \leq 0.1 w_{\min}.$$

The latter is impossible since even the $\alpha$-radius ball contains at least $3|S|w_{\min}/4$ points. This implies that $g(\alpha + s\nu) \geq 0.9$ and now again using (16), we see that $|Q'| \geq 0.8|S_i|$ (note that $Q'$ is found in step 4). Thus from Claim 2 (noting that the proof of works for any subset of $S_i$), we get that $\beta'$ is a fairly good approximation to $\sigma_{i,\max}$:

(23)                    $$2.5\sigma_{i,\max}^2 \geq \beta' \geq 0.16\sigma_{i,\max}^2.$$

From the definition of $s$ in step 3 of the algorithm, it follows that there is some $y \in S_i$ with $|x - y| \geq \alpha + (s-2)\nu$. So, from (14), we have $\alpha + (s-2)\nu \leq \sqrt{2}(R_i + 2t\sigma_{i,\max})$. So, we have

(24)    $$\alpha + s\nu + 3\sqrt{\beta'}\left(\log\frac{|S|}{\delta} + 1\right) \leq \sqrt{2}(R_i + 2t\sigma_{i,\max}) + 2\sigma_{i,\max}$$
$$+ 4\sigma_{i,\max}\left(\log\frac{|S|}{\delta} + 1\right).$$

Thus from (20), no point of $S_j$, $j \in P \setminus \{i\}$, is contained in $B(x, \alpha + s\nu + 3\sqrt{\beta'}(\log\frac{|S|}{\delta} + 1))$. So the set removed from $T$ in step 5 is a subset of $S_i$.

Finally, using $g(\alpha + s\nu) \geq 9/10$, and isoperimetry (Corollary 4), we see that

$$g\left(\alpha + s\nu + 3\sqrt{\beta'}\left(\log\frac{|S|}{\delta} + 1\right)\right) \geq 1 - \frac{\delta}{10|S|^2},$$

which implies that $\eta(x) \leq \alpha + s\nu + 3\sqrt{\beta'}(\log\frac{|S|}{\delta} + 1)$. Thus, by (13), all of $S_i$ is in $B(x, \alpha + s\nu + 3\sqrt{\beta'}(\log\frac{|S|}{\delta} + 1))$. This completes the inductive proof of correctness. $\square$

Now we prove a lemma that was used above when we estimated $\sigma_{\max}$.

LEMMA 12. *Suppose $F$ is a (general) Gaussian in $\Re^n$. If $L$ is a set of independent identically distributed samples, each distributed according* *to $F$, then with probability at least $1 - \frac{\delta}{10}$, we have [with $\varepsilon = 20n(\sqrt{\log n} + \sqrt{\log(1/\delta)})/\sqrt{|L|}$], every vector $w$ satisfies*

$$(25) \quad \begin{aligned} E_F(w \cdot (x - E_F(x))^2)(1 - \varepsilon) &\leq E_S(w \cdot (x - E_F(x))^2) \\ &\leq E_F(w \cdot (x - E_F(x))^2)(1 + \varepsilon), \end{aligned}$$

*where $E_S$ denotes the "sample mean"; that is, it stands for $\frac{1}{|L|} \sum_{x \in L} \cdot$*

PROOF. We may translate by $-E_F(x)$ and without loss of generality assume that $E_F(x)$ is the origin. Suppose $Q$ is the square root of the inverse of the variance–covariance matrix of $F$. We wish to prove, for all vectors $w$,

$$E_F((w \cdot x)^2)(1 - \varepsilon) \leq E_S((w \cdot x)^2) \leq E_F((w \cdot x)^2)(1 + \varepsilon).$$

Putting $Q^{-1}w = u$ (noting that $Q$ is nonsingular and symmetric), this is equivalent to saying, for all vectors $u$,

$$E_F((u \cdot (Qx))^2)(1 - \varepsilon) \leq E_S((u \cdot (Qx))^2) \leq E_F((u \cdot (Qx))^2)(1 + \varepsilon).$$

But $Qx$ is a random sample drawn according to the standard normal, so it suffices to prove the statement for the standard normal. To prove it for the standard normal, we proceed as follows. First, for each coordinate $i$, we have that $E_F(|x_i|^2) = 1$ and using properties of the standard one-dimensional normal density, for any real $s > 0$,

$$\mathrm{Prob}(|E_S(|x_i|^2) - 1| \leq s) \geq 1 - e^{-|L|s^2/4}.$$

Now consider a pair $i, j \in \{1, 2, \ldots, n\}$, where $i \neq j$. The random variable $x_i x_j$ has mean 0 and variance 1. $E_S(x_i x_j)$ is the average of $N$ i.i.d. samples (each not bounded, but we may use the properties of the normal density again) concentrated about its mean:

$$\mathrm{Prob}(E_S|x_i x_j| \leq s) \geq 1 - e^{-|L|s^2/100}.$$

Putting $s = 10 \frac{\sqrt{\log n}}{\sqrt{|L|}}$, we see that all these $O(n^2)$ upper bounds hold simultaneously with probability at least $1 - \delta/n^8$.

Thus we have that the "moment" of inertia matrix $M$ of $S$ whose $i, j$ entry is $E_S(x_i x_j)$ has entries between $1 - \frac{1}{2}\varepsilon$ and $1 + \frac{1}{2}\varepsilon$ on its diagonal and the sum of the absolute values of the entries in each row is at most $\varepsilon/2$. Thus by standard linear algebra (basically arguments based on the largest absolute value entry of any eigenvector), we have that the eigenvalues of $M$ are between $1 - \varepsilon$ and $1 + \varepsilon$, proving what we want. □

**4. Max-likelihood estimation.** Now we describe an algorithm for max-likelihood fit of a mixture of $k$ spherical Gaussians of equal radius to (possibly) unstructured data. First we derive a combinatorial characterization of the optimum solution in terms of the *k-median* (*sum of squares, Steiner version*) problem. In this problem, we are given $M$ points $x_1, x_2, \ldots, x_M \in \Re^n$ in $\Re^n$ and an integer $k$. The goal is to identify $k$ points $p_1, p_2, \ldots, p_k$ that minimize the function

$$(26) \qquad \sum_{i=1}^{M} |x - p_{c(j)}|^2,$$

where $p_{c(j)}$ is the point among $p_1, \ldots, p_k$ that is closest to $j$ and $|\cdot|$ denotes Euclidean distance.

THEOREM 13.    *The mixture of $k$ spherical Gaussians that minimizes the log-likelihood of the sample is exactly the solution to the above version of k-median.*

PROOF.    Recall the density function of a spherical Gaussian of variance $\sigma$ (and radius $\sigma\sqrt{n}$) is

$$\frac{1}{(2\pi\sigma)^{n/2}} \exp\left( -\frac{|x - p|^2}{2\sigma^2} \right).$$

Let $x_1, x_2, \ldots, x_M \in \Re^n$ be the points. Let $p_1, p_2, \ldots, p_k$ denote the centers of the Gaussians in the max-likelihood solution. For each data point $x_j$ let $p_{c(j)}$ denote the closest center. Then the mixing weights of the optimum mixture $w_1, w_2, \ldots, w_k$ are determined by considering, for each $i$, the fraction of points whose closest center is $p_i$.

The log-likelihood expression is obtained by adding terms for the individual points to obtain

$$-\left[ \text{Constant} + \frac{Mn}{2}\log\sigma + \sum_j \frac{|x_j - p_{c(j)}|^2}{2\sigma^2} \right].$$

The optimum value $\hat{\sigma}$ is obtained by differentiation,

$$(27) \qquad \hat{\sigma}^2 = \frac{2}{Mn} \sum_j |x_j - p_{c(j)}|^2,$$

which simplifies the log-likelihood expression to

$$\text{Constant} + \frac{Mn}{2}\log\hat{\sigma} + \frac{Mn}{4}.$$

Thus the goal is to minimize $\hat{\sigma}$, which from (27) involves minimizing the familiar objective function from the sum-of-squares version of the $k$-median problem.  $\square$

We indicate how to use known results about the $k$-median to provide a constant *additive factor* approximation to the log-likelihood. Charikar, Guha, Tardos and Shmoys [4] provide an $O(1)$ approximation to the $k$-median problem with sum-of-squares distance. They do not indicate if their algorithm works when the centers are not one of the sample points. However, the triangle inequality implies that picking the center to be one of the sample points changes the objective of the $k$-median problem by at most a factor 4. Hence we obtain a constant factor approximation to $\hat{\sigma}^2$, and hence an approximation to $\log(\hat{\sigma})$ that is correct within an $O(1)$ *additive* factor. More efficient algorithms for $k$-median are now known, so there is some hope that the observations of this section may lead to some practical learning algorithms.

**5. Conclusions.** Several open problems remain. The first concerns solving the classification problem for Gaussians with significant overlap. For example, consider mixtures of spherical Gaussians with pairwise intercenter distance only $O(\max\{\sigma_1, \sigma_2\})$. In this case, a constant fraction of their probability masses overlap, and the solution to the classification problem is not unique. Our algorithm does not work in this case, though a recent spectral technique of Vempala and Wang [25] does apply. (However, it does not apply to nonspherical Gaussians.)

The second problem concerns general Gaussians whose probability masses do not overlap much but which appear to coalesce under random projection. For example, consider a pair of concentric Gaussians that have the same axis orientation. (Of course, these axes are unknown and are not the same as the coordinate axes.) In $n-2$ axis directions their variance is $\sigma^2$, and in the other remaining two directions their variances are $1, \sigma$ and $\sigma, 1$, respectively. If $\sigma = \Omega(\log n)$, the difference in the last two coordinates is enough to differentiate their samples with probability $1 - 1/\text{poly}(n)$. But after projection to an $O(\log n)$-dimensional subspace, this difference disappears. Hence neither distance-based clustering nor projection-based clustering seems able to distinguish their samples.

The third open problem concerns max-likelihood estimation, which seems to involve combinatorial optimization with very bizarre objective criteria once we allow nonspherical Gaussians.

We suspect all the above open problems may prove difficult.

We note that Dasgupta (personal communication) has also suggested a variant of the classification problem in which the sample comes from a "noisy" Gaussian. Roughly speaking, the samples come from a mixture of sources, where each source is within distance $\varepsilon$ of a Gaussian. We can solve this problem in some cases for small values of $\varepsilon$, but that will be the subject of another paper. Broadly speaking, the problem is still open.

**Acknowledgments.** We thank Sanjoy Dasgupta for many helpful discussions, including drawing our attention to [21]. We also thank Laci Lovász and David McAllester for useful suggestions. We thank the anonymous referee for a quick, thorough critical reading and many important suggestions which helped improve our presentation.

## REFERENCES

[1] APPLEGATE, D. and KANNAN, R. (1991). Sampling and integration of near log-concave functions. In *Proceedings of the 23rd ACM Symposium on Theory of Computing* 156–163. ACM Press, New York.

[2] ARORA, S. and KANNAN, R. (2001). Learning mixtures of arbitrary Gaussians. In *Proceedings of the 33rd ACM Symposium on Theory of Computing* 247–257. ACM Press, New York.

[3] BOURGAIN, J. (1999). Random points in isotropic convex sets. In *Convex Geometric Analysis* 53–58. Cambridge Univ. Press. MR1665576

[4] CHARIKAR, M., GUHA, S., TARDOS, E. and SHMOYS, D. (1999). A constant-factor approximation algorithm for the $k$-median problem. In *Proceedings of the 31st ACM Symposium on Theory of Computing* 1–10. ACM Press, New York. MR1797458

[5] DASGUPTA, S. (1999). Learning mixtures of Gaussians. In *Proceedings of the 40th Annual IEEE Symposium on Foundations of Computer Science* 634–644. IEEE, New York. MR1917603

[6] DASGUPTA, S. and SCHULMAN, L. (2000). A two-round variant of EM for Gaussian mixtures. In *Proceedings of the 16th Annual Conference on Uncertainty in Artificial Inteligence* 152–159. Morgan Kaufmann, San Francisco, CA.

[7] DEMPSTER, A. P., LAIRD, N. M. and RUBIN, D. B. (1977). Maximum likelihood from incomplete data via the EM algorithm. *J. Roy. Statist. Soc. Ser. B* **39** 1–38. MR501537

[8] DRINEAS, P., FRIEZE, A., KANNAN, R., VINAY, V. and VEMPALA, S. (1999). Clustering in large graphs and matrices. In *Proceedings of the ACM–SIAM Symposium on Discrete Algorithms* 291–299. ACM Press, New York. MR1739957

[9] DYER, M. E., FRIEZE, A. and KANNAN, R. (1991). A random polynomial time algorithm for approximating the volume of convex bodies. *Journal of the ACM* **38** 1–17. MR1095916

[10] FREUND, Y. and MANSOUR, Y. (1999). Estimating a mixture of two product distributions. In *Proceedings of the 12th ACM Conference on Computational Learning Theory* 53–62. ACM Press, New York. MR1811601

[11] HOEFFDING, W. J. (1963). Probability inequalities for sums of bounded random variables. *J. Amer. Statist. Assoc.* **58** 13–30. MR144363

[12] HUBER, P. J. (1985). Projection pursuit. *Ann. Statist.* **13** 435–475. MR790553

[13] JOHNSON, W. B. and LINDENSTRAUSS, J. (1984). Extensions of Lipshitz mapping into Hilbert space. *Contemp. Math.* **26** 189–206. MR737400

[14] KANNAN, R. and LI, G. (1996). Sampling according to the multivariate normal density. In *Proceedings of the 37th Annual IEEE Symposium on Foundations of Computer Science* 204–213. IEEE, New York. MR1450618

[15] LINDSAY, B. (1995). *Mixture Models: Theory, Geometry, and Applications.* Amer. Statist. Assoc., Alexandria, VA.

[16] Lovász, L. and Simonovits, M. (1990). Mixing rate of Markov chains, an isoperimetric inequality, and computing the volume. In *Proceedings of the 31st Annual IEEE Symposium on Foundations of Computer Science* 346–355. IEEE, New York. MR1150706

[17] Lovász, L. and Simonovits, M. (1993). Random walks in a convex body and an improved volume algorithm. *Random Structures Algorithms* **4** 359–412. MR1238906

[18] Megiddo, N. (1990). On the complexity of some geometric problems in unbounded dimension. *J. Symbolic Comput.* **10** 327–334. MR1086309

[19] Prékopa, A. (1971). Logarithmic concave measures with applications to stochastic programming. *Acta Sci. Math.* (*Szeged*) **32** 301–316. MR315079

[20] Prékopa, A. (1973). On logarithmic concave measures and functions. *Acta Sci. Math.* (*Szeged*) **34** 335–343. MR404557

[21] Redner, R. A. and Walker, H. F. (1984). Mixture densities, maximum likelihood and the EM algorithm. *SIAM Rev.* **26** 195–239. MR738930

[22] Rudelson, M. (1999). Random vectors in the isotropic position. *J. Funct. Anal.* **164** 60–72. MR1694526

[23] Titterington, D. M., Smith, A. F. M. and Makov, U. E. (1985). *Statistical Analysis of Finite Mixture Distributions.* Wiley, New York. MR838090

[24] Vapnik, V. N. and Chervonenkis, A. Ya. (1971). Theory of uniform convergence of frequencies of events to their probabilities and problems of search to an optimal solution from empirical data. *Automat. Remote Control* **32** 207–217. MR301855

[25] Vempala, S. and Wang, G. (2002). A spectral algorithm for learning mixture models. In *Proceedings of the 43rd Annual IEEE Symposium on Foundations of Computer Science* 113–124. IEEE, New York.

Department of Computer Science
Princeton University
Princeton, New Jersey
USA
e-mail: arora@cs.princeton.edu

Department of Computer Science
Yale University
51 Prospect Street
New Haven, Connecticut 06517
USA
e-mail: kannan@cs.yale.edu